# Uso de gestos –como recurso-mediador– por un profesor de bachillerato para enfrentar un desafío didáctico no previsto por él

Ulises Alfonso Salinas-Hernández, Luc Trouche

**Resumen**

En este artículo se reporta cómo un profesor de bachillerato (niveles 10-12) usó gestos al enfrentarse a una situación en el aula no prevista por él. El uso de gestos ocurrió durante la interacción con tres estudiantes en torno al análisis del movimiento de un objeto en caída libre. Se llevó a cabo un análisis cualitativo de los datos a partir de un marco conceptual que coordina elementos de tres aproximaciones teóricas: aproximación documental de lo didáctico, la teoría de actividad y la epistemología histórica. La recopilación de datos se llevó a cabo mediante la videograbación de la interacción profesor-estudiantes al momento de tratar de dar significado al movimiento de un objeto en caída libre. Consideramos que el análisis de resultados muestra que la coordinación de gestos es, por una parte, un recurso semiótico que debe considerarse relevante en el sistema de recursos de los profesores; y por otra parte, representan, en sí mismo, una manera de visualizar la estructura cognitiva que guía las acciones del profesor.
**Palabras clave**: Recursos; gestos; artefactos; teoría de la actividad; mediación

**Abstract**

In this article, we report how a high school teacher (levels 10-12) used gestures when facing a classroom situation not foreseen by him. The use of gestures took place during the interaction with three students around the analysis of the movement of an object in free fall. A qualitative analysis of the data was carried out from a conceptual framework that coordinates elements of three theoretical approaches: documentational approach to didactics, activity theory and historical epistemology. The data collection was carried out by videotaping the teacher-student interaction while the meaning production process about the movement of an object in free fall took place. We consider that the analysis of results shows that the coordination of gestures is, on the one hand, a semiotic resource that should be considered relevant in the teachers' resources system; and on the other hand, represents in itself, a way of visualizing the cognitive structure that guides the actions of the teacher.

**Résumé**

Dans cet article, nous montrons comment un professeur de lycée (grades 10-12) utilise des gestes pour faire face à une classe de situations inattendue pour lui. Ces gestes prennent place à l'occasion d'un échange avec trois étudiants pendant la description de la trajectoire d'un objet en chute libre. L'analyse qualitative des données est conduite à partir d'un cadre théorique qui coordonne trois approches: l'approche documentaire du didactique, la théorie de l'activité, et l'épistémologie historique. Les données sont recueillies à partir de vidéos capturées au cours du processus de construction de signification de la trajectoire d'un objet en chute libre. Nous considérons que l'analyse des résultats met en évidence que la coordination des gestes est, à la fois, une ressource sémiotique partie prenante du système de ressources des enseignants, et un moyen de visualiser la structure cognitive qui guide les actions du professeur.
**Mots clés:** Ressources, gestes, artefacts, activité, médiation

## 1. Introducción y problema de investigación

En el terreno didáctico, trabajos sobre la práctica del profesor de matemáticas coinciden en la necesidad de poner atención en el *uso* de *recursos*; entender qué son y la manera en que funcionan como extensiones del profesor durante su práctica docente (Adler, 2000). Además, esta autora señala que los *recursos* no se restringen a objetos materiales. Ella los categoriza [a los *recursos*] en humanos (profesores, padres de familia, conocimiento del profesor, entre otros), materiales (libros de texto, calculadoras y objetos matemáticos como plano cartesiano, entre otros) y socio-culturales (lenguaje). Guzmán y Kieran (2013) señalan que la manera en que los *recursos* apoyan o no a los profesores en sus esfuerzos de resolver problemas en clase, claramente tienen un impacto en la experiencia de resolver problemas de los estudiantes. Se retoman las investigaciones sobre el *uso* de *recursos* en profesores de matemáticas para vincularlas con la necesidad de analizar el *uso* de *recursos* en la práctica de los profesores de física. Se parte del hecho de la estrecha relación de la física y las matemáticas como disciplinas científicas, tanto en el plano epistemológico como en el plano educativo. En el plano epistemológico, respecto de la relación entre física y matemáticas, como disciplinas científicas, Piaget señala:

> Desde el primer contacto con la epistemología física nos encontramos, pues, en presencia de la dificultad sumamente instructiva de delimitar los campos entre la física y la matemática: o reducimos los dos a uno solo, o nos empeñamos en distinguirlos, pero sin alcanzar una frontera estática. (Piaget, 1979, p. 9).

En el ámbito escolar del campo de la matemática educativa, las líneas de investigación referentes al papel de las matemáticas en las clases de física, o viceversa, son diversas (e.g. Karam, 2015; Kragh, 2015; Kjeldsen & Lützen, 2015). Por su parte, dentro del interés de la presente investigación sobre el *uso* de gestos por parte del profesor, se sigue la perspectiva epistemológica planteada por Radford (2009), quien señala, por una parte, que el pensamiento matemático no sólo está mediado por símbolos escritos, sino también por las acciones, gestos y otros tipos de señales. Y, por otra parte, que el pensamiento se produce también a través de una sofisticada coordinación semiótica de la voz, el cuerpo, los gestos, los símbolos y las herramientas[1]. Así, empleando el término *artefacto* para incluir a los *recursos* didácticos, Radford (2012) señala la importancia de investigar su *uso* y comprender su influencia en los procesos de enseñanza y de aprendizaje. De manera particular,

---

[1] Esta perspectiva se asocia con el paradigma multimodal, el cual cuestiona la hegemonía de la representación y la comunicación a través del código escrito para el aprendizaje y destaca el papel de otros sistemas semióticos para construir significados (Haquin, 2012).

Radford (2006) expone que –a diferencia de las aproximaciones racionalistas en las cuales el pensamiento corresponde a una actividad mental en el que la mente no necesita la asistencia de los sentidos ni de la experiencia para alcanzar las verdades matemáticas– el pensamiento es una práctica social en el sentido de Wartofsky (1979). En palabras de Radford (2006): "[E]l pensamiento es considerado una reflexión mediatizada del mundo de acuerdo con la forma o modo de la actividad de los individuos" (p. 107).

Así, a partir, por un lado, de la importancia en analizar el *uso* de *recursos* por los profesores en el salón de clase cuando resuelven problemas; y por otro lado, siguiendo una concepción no mentalista del pensamiento, en este artículo pretendemos responder las siguientes preguntas de investigación: (1) ¿cómo son utilizados los gestos por parte del profesor, como *recursos* y como modos de representación [*artefactos*] para que un estudiante interprete el movimiento de un objeto en caída libre en función del *sistema de referencia*? (2) ¿qué se puede inferir sobre el conocimiento del profesor –en torno a la descripción de la caída libre de objetos– a partir de sus gestos?

El artículo se divide en seis secciones. Después de haber presentado la introducción y el problema de investigación (§ 1), se define el marco conceptual en el que se integra la teoría de la actividad (Engeström, 2001), el *uso* de *recursos* (Gueudet & Trouche, 2009, 2012) y la producción de *artefactos* (Wartofsky, 1979) (§ 2). Se continúa con un análisis *a priori* del contenido matemático-físico que se presenta en el artículo (§ 3); seguido por la metodología (§ 4) y el análisis de resultados (§ 5). Las conclusiones y reflexiones finales se presentan en la última sección del artículo (§ 6).

## 2. Marco conceptual

En esta sección se presenta el encuadre teórico del artículo, el cual tiene como objetivo atender tanto la dimensión didáctica, como la dimensión epistemológica de las preguntas de investigación (§ 1). La dimensión didáctica, que se atiende a través de la *Aproximación documental de lo didáctico* (ADD) (Gueudet & Trouche, 2009, 2012), se asocia con el *uso* de gestos por parte del profesor para afrontar la situación pedagógica que se le presenta –a partir de sus conocimientos y de la manera en que tiene integrado el *uso* de gestos en su *sistema de recursos*–. Mientras que la dimensión epistemológica, que se atiende con la noción de *artefactos* de Wartofsky (1979), se relaciona con la manera en que el conocimiento del profesor es representado a través de los gestos. Lo que conlleva, a su vez, en analizar la componente epistemológica de los gestos. Las aproximaciones antes señaladas se coordinan a través de la teoría de la actividad (TA) (Engeström, 2001). En el caso particular de las investigaciones sobre la práctica del profesor, trabajar conjuntamente con diferentes aproximaciones teóricas ayuda a ser consciente tanto de los aspectos que puedan faltar en tales aproximaciones, como también de los beneficios que puedan aportar cada una con el objetivo de tener una mejor comprensión de la práctica del profesor (Trouche, Gitirana, Miyakawa, Pepin, & Wang, Online first).

Gueudet y Trouche (2009; 2012) proponen un enfoque teórico cercano, en la conceptualización de *recursos*, de la propuesta de Adler (2000). Por ejemplo, no concebir a los *recursos* sólo como los provenientes de objetos materiales, sino a todos aquellos que intervienen en la comprensión y resolución de problemas. En su propuesta, estos autores hacen la diferencia entre *recursos* y documentos. Así, los documentos son desarrollados a través de lo que denominan *génesis documental*. El *trabajo documental* es el núcleo de la actividad de los profesores y de su desarrollo profesional. Gueudet y Trouche (2009; 2012) utilizan el término *recurso* para dar énfasis a la variedad de *artefactos* que consideran y en donde a su vez un *artefacto* (físico o psicológico) es un medio cultural y social provisto por la actividad humana (e.g., computadoras y lenguaje); producidos con propósitos específicos (e.g., resolver problemas).

En la *génesis documental*, los documentos son creados a partir de un proceso en el cual los profesores construyen esquemas de utilización de los *recursos* para situaciones dentro de una variedad de contextos, proceso que se ejemplifica por la ecuación: *Documento = recursos + esquemas* de utilización (Gueudet & Trouche, 2009, p. 205). Los mismos autores mencionan que los esquemas de utilización supone una parte observable y otra invisible. Los *usos*–las reglas particulares de acción– corresponden a la parte observable del esquema, que ocurre durante las acciones del profesor en el transcurso de la actividad. Mientras que las invariantes operatorias corresponden a la estructura cognitiva que guía las acciones del profesor. Así, los esquemas sólo son observables a través de las acciones [*usos*] que lleva a cabo el sujeto al trabajar con los *recursos*. En palabras de Gueudet y Trouche (2009): "Entonces, el investigador puede tratar de inferir invariantes operatorias de los *usos* (p. 209; traducción libre). Se tiene así una segunda relación: *Documento=Recursos+Usos+Invariantes operatorias*.

Durante la actividad del sujeto con el *uso* de los *artefactos* se dan dos procesos: instrumentación e instrumentalización; el primero tiene que ver con la influencia del *artefacto* en las acciones [actividad] del sujeto. Mientras que la instrumentalización ocurre cuando el sujeto se apropia del *artefacto* y determina la manera en que se usa. Es importante resaltar que es durante la instrumentalización que el sujeto [e.g., el profesor] adapta y modifica los *artefactos* [*recursos*] de acuerdo a la variedad de situaciones que se le presentan. Así, es durante este proceso que se desarrolla la componente creativa-didáctica del profesor para enfrentar un problema.

En síntesis, al estudiar el desarrollo de la *génesis documental* de los profesores se obtienen evidencias de la manera en que ellos articulan los diferentes documentos. Así, para dar cuenta de la variedad y del diseño de documentos y la manera en que los profesores los articulan en una variedad de situaciones, los documentos se estructuran en un *sistema de documentación*, en donde el *sistema de recursos* del profesor constituye la parte del "*recurso*" del *sistema de documentación* –sin considerar la parte del esquema del documento– (Gueudet & Trouche, 2012).

Wartofsky (1979) desarrolló su aproximación teórica, la cual denominó Epistemología histórica (EH), con la intención de mostrar cómo los modos más evolucionados de representación que ha logrado el ser humano (e.g., las teorías científicas) tienen su génesis en los modos de representación que surgen simultáneamente con nuestra práctica productiva, social y lingüística primaria. Para él la característica fundamental de la práctica cognitiva humana es la habilidad para crear representaciones; además los seres humanos crean los medios de su propia cognición: los *artefactos*. El *artefacto* para Wartofsky (1979) es tanto un medio cognitivo como un modo de representación. Es decir, representan el modo de actividad en el que fueron producidos y son además un medio de transmitir y lograr conocimiento. Así, además de considerar al lenguaje, Wartofsky (1979) introduce las formas de organización e interacción social, las técnicas de producción, y la adquisición de habilidades, como *artefactos*. De manera que, al producir *artefactos* [medios cognitivos] para su *uso*, se producen representaciones. Las teorías científicas, los sistemas lógicos y las matemáticas (además de las formas de representación en literatura y el arte) son los modos más evolucionados de representación que ha logrado el ser humano. Y es a través de las representaciones como el ser humano logra el conocimiento. De manera que el *uso* y la creación de *artefactos* es una forma de acción [*praxis*] distintivamente humana.

Se observa que tanto en la ADD como en la HE un elemento esencial es conceptualizar las prácticas humanas como una práctica social orientada por objetivos. De manera que se combinan ambas aproximaciones teóricas a través de la teoría de la actividad (TA) en la forma desarrollada por Engeström (2001). Lo cual nos permite integrarlas en marco conceptual común. Engeström (2001) extendió el planteamiento de Vygotsky sobre la mediación a través de los signos y herramientas y el de Leont'ev (1978) sobre la labor para presentar, como unidad de análisis, al menos dos *sistemas de actividad* interactuando entre sí.

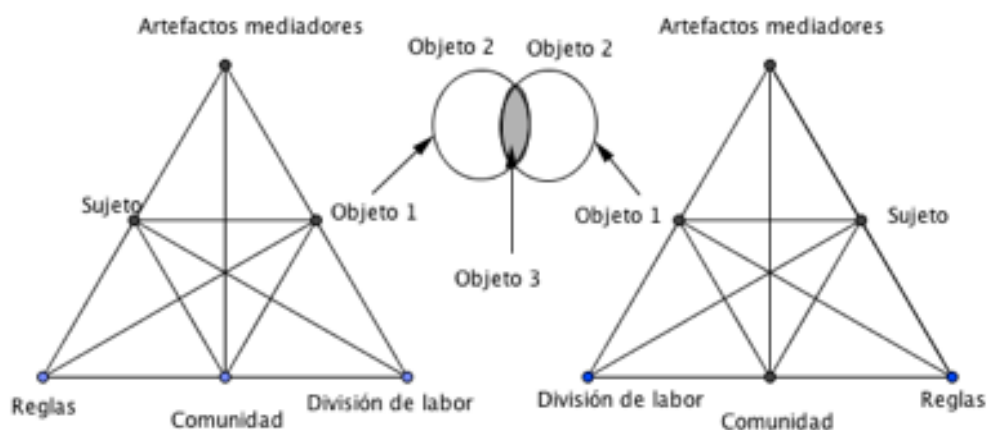

**Figura 2.1. Representación de Engeström sobre el *sistema de actividad*, basado en la mediación y en el trabajo de Leont'ev.**

La Figura 2.1 presenta los principios centrales que sintetizan la aproximación: (1) Un *sistema de actividad* colectiva, mediado por *artefactos* y orientado a objetos es la unidad de análisis principal. (2) Hay una variedad de voces en los *sistemas de actividad*. Así, un sistema de actividad es siempre una comunidad con diferentes puntos de vista, tradiciones e intereses. Además, la división de labor en una actividad crea diferentes posiciones para los participantes. (3) Existe una historicidad en los *sistemas de actividad*. Estos toman forma y son transformados durante grandes periodos. (4) Las contradicciones tienen un rol importante en el cambio y desarrollo de los *sistemas de actividad*. (5) Existe la posibilidad de transformaciones expansivas en los *sistemas de actividad*. Ellos se mueven a través de ciclos relativamente largos de transformaciones cualitativas.

A partir de lo presentado en el marco conceptual, se señala que la práctica docente es, además de una práctica pedagógica guiada por el *trabajo documental* del profesor y el *uso* de *recursos*, una actividad social de significación llevada a cabo por medio de signos y *artefactos* [e.g., gestos] que alteran de manera fundamental la forma como pensamos y actuamos. De manera que el análisis de datos va encaminado en determinar las componentes didáctica y epistemológica de los gestos dentro de los *sistemas de actividad* que se desarrollan durante la discusión del profesor con un estudiante; a partir, por un lado, de su *uso* como *recurso*; y por el otro lado, como signos y *artefactos* que representan el modo de actividad del profesor, esto es, su conocimiento.

## 3. Metodología

Este estudio es de carácter cualitativo y se llevó a cabo en un laboratorio de física de un bachillerato (nivel 10-12) de la ciudad de México. Participaron: el profesor, quien imparte la asignatura de física desde hace más de 30 años y 11 estudiantes (16-18 años). Los datos analizados en la presente investigación fueron recabados durante la recopilación de datos de otra investigación (Salinas-Hernández & Miranda, 2018), en la cual no fueron analizados e incorporados. En la investigación de Salinas-Hernández y Miranda (2018), se diseñaron cinco tareas para que fueran llevadas a cabo por los estudiantes en equipos de 3 y 4 integrantes. Las tareas fueron diseñadas por los investigadores y mostradas al profesor antes de ser implementadas, con el objetivo de que el profesor tuviera conocimiento sobre lo que estaban por realizar sus estudiantes. El papel del profesor consistió en supervisar el trabajo de los estudiantes mientras resolvían las tareas y aclarar algunas dudas que pudieran surgir.

En conjunto, las tareas tuvieron el propósito de analizar en los estudiantes el proceso de interpretación de gráficas cartesianas obtenidas por ellos a partir de un software [sensor virtual de movimiento], y que estaban relacionadas con un experimento de la caída de una pelota por un plano inclinado. Así, las tareas se dividieron en dos partes: las primeras dos corresponden a la implementación del experimento –de la caída de una pelota por un plano inclinado– por parte de los estudiantes; mientras que las tareas 3, 4 y 5 correspondieron a preguntas (en hojas de trabajo) para indagar la comprensión de los estudiantes entorno al fenómeno físico de la caída de un móvil por un plano inclinado y la relación con las gráficas cartesianas generadas por dicho movimiento. Durante el periodo de recolección de datos todos los participantes fueron videograbados con dos cámaras, éstas se enfocaron principalmente a los estudiantes mientras llevaban a caba cada una de las tareas.

Para los propósitos de este artículo se presenta el análisis de un extracto de video de 10:38 minutos de duración, el cual fue tomado mientras un equipo de tres estudiantes resolvía la Tarea 4 (Figura 3.1). En este extracto de video surge una pregunta (la cual no aparecía en ninguna de las tareas) en un estudiante [E1] relacionada con el concepto de aceleración negativa. El análisis de datos está enfocado en la manera en que el profesor [Prof.] trata de resolver la duda de los estudiantes a través del *uso* de gestos y establece un diálogo con uno de los estudiantes de ese equipo a quien se denominó Pedro.

Se deja fijo un plano inclinado (Figura 10) y se baja la altura del otro plano inclinado hasta una altura de 60 *cm* (Figura 11).

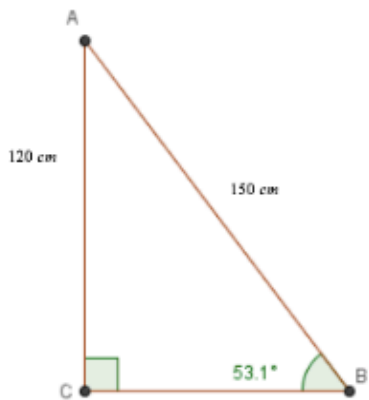

Figura 10: Plano inclinado P1

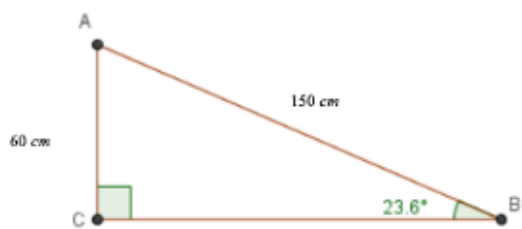

Figura 11: Plano inclinado P4

PREGUNTAS:

1. ¿Desde qué plano inclinado crees que llegue primero la pelota al suelo? Argumenta tu respuesta.

**Figura 3.1. Tarea (4) que se encontraban realizando los estudiantes que aparecen en el extracto analizado.**

## 4. Contenido conceptual-didáctico de la investigación y análisis *a priori*

La presente sección se divide en dos partes. En la primera parte se presenta el contenido conceptual-didáctico (físico-matemático) de la investigación. En la segunda, se lleva a cabo un análisis *a priori* de dicho contenido considerando 3 niveles: (1) epistemológico, (2) didáctico-curricular y (3) relacionado con el *uso* de *recursos*.

### 4.1. Contenido conceptual-didáctico de la investigación

Para dar respuesta a las preguntas de investigación de este artículo (§ 1) es necesario dar cuenta de los dos conceptos principales que están en juego en el proceso de significación del movimiento de caída libre que ocurre durante el diálogo profesor-estudiante. Estos conceptos son el *sistema de coordenadas cartesianas* y el concepto físico de *sistema de referencia*.

#### 4.1.1. *Sistema de coordenadas cartesianas*

El *sistema de coordenadas cartesianas* (SCC) (Figura 4.1.) se forma por dos ejes ortogonales denominados habitualmente: eje X (eje horizontal o de las abscisas) y eje Y (eje vertical o de las ordenadas). Esta representación es una manera de identificar la posición de cualquier punto «par ordenado con coordenadas (x,y)» con respecto del origen de coordenadas (punto de intersección de los ejes). Se divide en cuatro cuadrantes y los valores (positivos o negativos) tanto de las abscisas como de las ordenadas se determinan de acuerdo con su posición respecto del origen de coordenadas. Así, en el eje X las abscisas toman valores positivos a la derecha del origen de coordenadas y negativos a la izquierda del mismo. Mientras que en el eje Y, las ordenadas toman valores positivos hacia arriba del origen de coordenadas y negativos hacia abajo del mismo punto.

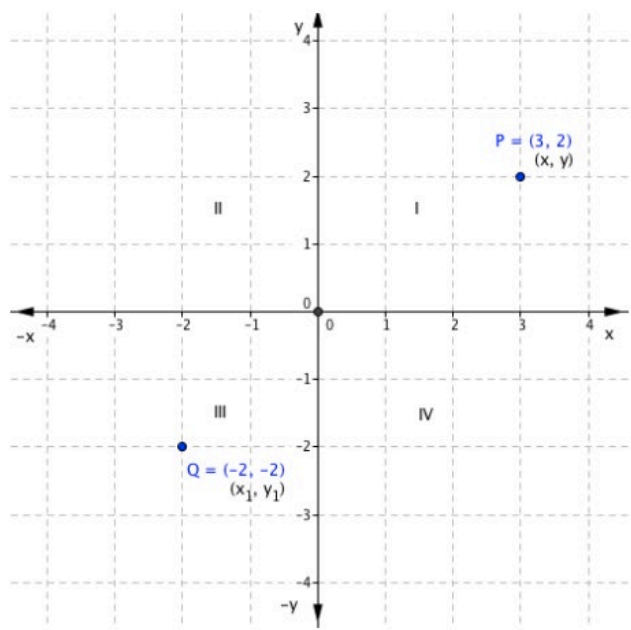

**Figura 4.1. Representación del *sistema de coordenadas cartesianas*.**

## 4.1.2. *Sistema de referencia*

En física, y en particular en mecánica newtoniana, se entiende como *sistema de referencia* (SR) o *marco de referencia* al conjunto de convenciones usadas por un observador (respecto del cual se realizan las mediciones) para poder medir magnitudes físicas de un sistema físico (e.g., posición, velocidad, aceleración, distancia) y en el cual, a su vez, se cumplen las tres leyes de Newton. En un texto clásico, para los primeros ciclos universitarios en México, se puntualiza que los observadores en los diferentes *sistemas de referencia* pueden obtener distintos valores numéricos de las cantidades físicas medidas, pero las relaciones entre las cantidades [las leyes de la física], son las mismas para todos los observadores (Resnick & Halliday, 1970). Sobre su importancia, los mismos autores dicen: "Por consiguiente, es importante que el estudiante siempre se dé cuenta de cuál es su marco de referencia en un determinado problema." (p. 30).

## 4.2. Análisis *a priori* de los conceptos: *sistema de coordenadas cartesiano* y *sistema de referencia*

Se lleva a cabo ahora el análisis *a priori* de los dos conceptos (SCC y SR) considerando 3 niveles: (1) epistemológico, (2) didáctico-curricular y (3) relacionado con el *uso* de *recursos*. El nivel epistemológico considera lo expuesto anteriormente, por un lado, sobre la estrecha relación entre física y matemáticas como prácticas científicas (§ 1), y por otro lado, sobre el contenido conceptual-didáctico de la investigación (§ 4.1), para señalar que, en particular, existe una estrecha relación epistemológica entre el significado físico del SR y el significado matemático del SCC. Incluso, en el terreno epistemológico-didáctico, los SR suelen ser representados gráficamente a partir del SCC, y las gráficas que se generan para interpretar el movimiento de objetos, suelen ser denominadas: gráficas cartesianas. Sin embargo, es importante tener en cuenta que la representación gráfica de un SR no es un SCC. Los significados de ambos conceptos fueron ya señalados (§ 4.1.1 y § 4.1.2). Es aquí en donde se puede presentar un obstáculo epistemológico al asociar el significado del *sistema de referencia* con la representación del *sistema de coordenadas cartesianas* (Figura 4.1). El nivel didáctico-curricular, toma en cuenta que el SCC y el SR forman parte de los contenidos curriculares esenciales para ser abordados por parte de los profesores en las asignaturas de matemáticas y física, respectivamente (Figura 4.2).

| LÍNEAS TEMÁTICAS | TEMÁTICA |
|---|---|
| **Eje 4:**<br>**Geometría Analítica.**<br>Sistema de coordenadas. Plano Cartesiano. Estudio analítico de problemas de corte euclidiano y de lugares geométricos.<br><br>**Eje 5: Funciones y Plano Cartesiano.**<br>Concepto de función y sus elementos. Diversos tipos de variación, estudio de sus comportamientos. Relación parámetro- gráfica- variación. | **1. Primera Ley de Newton**<br>• Inercia, sistema de referencia y reposo.<br><br>• Interacciones y fuerzas, aspecto cualitativo.<br><br>• Fuerza resultante cero, (vectores desde un punto de vista operativo, diferencia entre vector y escalar), 1ª Ley de Newton y Movimiento Rectilíneo Uniforme. |

**Figura 4.2. Contenidos curriculares del bachillerato en el que trabaja el profesor de este estudio: SCC en matemáticas (izquierda) y SR en física (derecha).**

El nivel sobre el *uso* de *recursos* asocia el SCC y el SR con los posibles *recursos* para ser movilizados por parte del profesor para enfrentar el desafío didáctico, esto es, la pregunta de un estudiante sobre la aceleración negativa (§ 3). Así, se consideran como posibles *recursos* el *uso*: del pizarrón, en el cual el profesor pueda hacer gráficos (e.g. usar el propio SCC); usar otro ejemplo que se relacione con la pregunta del estudiante y gestos. El *uso* de *recursos* está dirigido a relacionar el significado de la aceleración negativa con el *uso*, a su vez, de un SR (elegido arbitrariamente).

## 5. Análisis de resultados

A continuación se presenta el análisis del extracto del video de 10:38 minutos de duración. Se analiza el diálogo que ocurre entre el profesor y los estudiantes a partir de la interacción de 3 *sistemas de actividad*. Cada *sistema de actividad* está integrado por los mismos sujetos (el equipo de 3 estudiantes y el profesor), pero orientado, cada uno, a objetivos particulares distintos. El primer *sistema de actividad* (SA1, analizado en Salinas-Hernández & Miranda, 2018) está orientado a resolver la Tarea 4 por parte de los estudiantes. Este objetivo surge a partir de los intereses de la investigación de Salinas-Hernández y Miranda (2018). El segundo *sistema de actividad* (SA2) surge a partir de la pregunta de E1, es decir, el objetivo (determinado por los estudiantes) es saber la diferencia entre aceleración negativa y desaceleración. El tercer *sistema de actividad* (SA3) se determina a partir del interés del profesor por hacer, a partir de la pregunta de E1, que los estudiantes interpreten el movimiento de caída libre de un objeto en términos de un SR elegido arbitrariamente. Sin embargo existe un objetivo común (e.g., Objeto 3 de la Figura 2.1), el cual está dado por el papel que juega el SR en la interpretación del movimiento de objetos (y por tanto se dé significado a conceptos de posición, velocidad y aceleración).

De esta manera, el análisis está enfocado en determinar de qué manera el profesor relaciona los tres *sistemas de actividad*, a través del *uso* de gestos. Este análisis, a su vez, permite determinar cómo está integrado el *uso* de gestos y los conceptos de SCC y SR en el sistema de *recursos* del profesor. Así, el video se divide, para su análisis, en tres momentos que están determinados por el desarrollo de SA2 (§ 5.1), de SA3 (§ 5.2) y la integración (por parte del profesor) de los tres *sistemas de actividad* (5.3).

### 5.1. Inicio de la discusión: *segundo sistema de actividad*

El SA2 tiene lugar cuando a los estudiantes, quienes se encontraban en el SA1 atendiendo los objetivos de otros sujetos (Salinas-Hernández & Miranda, 2018), les surge un interés particular (objetivo del SA2): conocer el significado de la aceleración negativa. Y es cuando inicia el diálogo con el profesor.

L1     E1: ¿Cuál sería un ejemplo de aceleración negativa, profe?

L2     Pedro: Es que desaceleración es dejar de acelerar, o sea, que un móvil se vaya deteniendo poco a poco.

L3     E1: ¿Cuál sería un ejemplo de aceleración negativa? [*Dirigiéndose al profesor*].

L4     Prof.: Nada más se le llama desaceleración hasta el momento en que llegas a cero.

L5     Pedro: ¡Ah! Ok.

L6     Prof.: ¿Si?, estás desacelerando hasta llegar a cero. Que está asociado con la aceleración negativa. Nada más que la aceleración negativa no para, no cesa ahí. Continúa, ¿no? Por ejemplo, cuando tú avientas un objeto, desacelera.

L7     E1: Y ya va a tener una aceleración negativa, si pasas del cero.

L8   Prof.: Siempre tiene una aceleración negativa, siempre tiene una aceleración negativa nada más que la desaceleración tú terminas cuando se para [*se detiene*]. Y aquí sigue actuando de modo que ahora [*es interrumpido por E1 que dice: "nos da negativo"*]. Entonces, va el cuerpo para arriba, pero nuestra aceleración, negativa. Entonces, frena [*cuando llega a su punto más alto*], pero como la aceleración continúa, baja [*se refiere a que el cuerpo cambia de dirección, respecto de la inicial*].

El SA2 se desarrolla entre L1 y L8. Da inicio cuando E1 plantea una pregunta al profesor relativa a lo que es la aceleración negativa (L1). Aquí, el profesor retoma en L4 y L6 lo dicho por P. [Pedro] previamente (L2). Él [Pedro] trató de relacionar la aceleración negativa [concepto] con un fenómeno físico conocido y producto de su experiencia; sin embargo, el profesor no es claro en su explicación de la relación que existe entre la desaceleración y la aceleración negativa, no obstante que dice: "Que está asociado con" (L6). La relación que hace el profesor es la siguiente: cuando el profesor dice en L6 "estás desacelerando hasta llegar a cero" hace referencia al momento en que el móvil se detiene, es decir, al fenómeno; y al decir "Nada más que la aceleración negativa no para, no cesa ahí. Continua…" hace referencia a la aceleración negativa [concepto] que se obtiene a partir de un SR.

El profesor no deja claro cómo un concepto (como el de aceleración negativa) se puede emplear para analizar un fenómeno (desaceleración). En su discurso, el profesor argumenta (implícitamente) que la dirección del objeto ocurre hacia donde se ha determinado (a partir de un SR) que la aceleración tiene signo negativo. Pero pueden ocurrir aceleraciones con signo positivo (aceleraciones positivas) y que aún así un móvil está desacelerando (frenando), dependiendo del SR elegido por el observador. De ahí la importancia de la elección de un SR adecuado para analizar un determinado movimiento de objetos. En un libro de texto (Giancoli, 2006), que se utiliza como apoyo durante el desarrollo de los cursos de este nivel escolar, se manifiesta el cuidado que se debe de poner al no concebir que desaceleración significa necesariamente que la aceleración sea negativa. "Cuando un objeto frena, a veces se dice que está desacelerando. Pero hay que ser cautelosos: desaceleración no significa necesariamente que la aceleración sea negativa." (p. 25). Más bien, significa que la magnitud de la velocidad disminuye. Posteriormente, el profesor da un ejemplo y dice: "cuando tú avientas un objeto, desacelera". Aquí, el profesor no explicita si hace referencia al tiro vertical, o bien se infiere que se trata de dicho movimiento cuando dice: "Entonces va el cuerpo para arriba […]. Entonces, frena, pero como la aceleración continúa, baja". El ejemplo del profesor corresponde ahora a un movimiento que tiene un comportamiento en dos direcciones (cuando sube y cuando baja), pero en donde la orientación del SR que indica aceleración negativa hacia abajo es constante. Tampoco vincula la desaceleración con los cambios de velocidad del objeto. A partir de este momento, el profesor tiene como objetivo que los estudiantes comprendan cómo se establece el signo (positivo o negativo) de las cantidades físicas (e.g., aceleración), para lo cual hará *uso*, en particular, de dos *recursos*: el concepto de *sistema de referencia* y gestos; los cuales hasta este momento no han sido usados por el profesor. Para hacer *uso* de estos *recursos* (SR y gestos) el profesor retoma el experimento realizado por los estudiantes del plano inclinado (SA1) pero con otro objetivo, que da lugar al SA3.

### 5.2. Cambio hacia otro *sistema de actividad*: tercer *sistema de actividad*

El SA3 tiene lugar de L9 a L20, cuando el profesor orienta la discusión hacia otro sistema de actividad (SA3).

L9   Prof.: En el plano inclinado [*Refiriéndose al experimento que realizaron los estudiantes*], ¿para dónde tomaron lo positivo? [*Pedro pregunta: "¿cómo?"*], ¿para dónde tomaron lo positivo en el plano inclinado? [*Los estudiantes se ponen nerviosos*] Sí, pusieron positivo y negativo, ¿no? [*Pedro responde: "sí"*] ¿Para dónde estaba lo positivo?

L10  Pedro: Este…positivo, bueno, yo lo que logré entender, fue que positivo se tomó desde, bueno desde la parte, ¿cómo se dice?, desde la longitud más arriba, bueno los centímetros que había de la distancia del suelo hasta los ciento veinte centímetros que nos venía en la práctica. Y desde ese punto hasta el suelo, bueno a donde terminaba el riel [*refiriéndose al plano inclinado*] era, este, el plano.

En L9, evidencia el cambio entre *sistemas de actividad*. Al ponerse nerviosos (lo cual se observó en el video) los estudiantes dan muestra que la pregunta del profesor (L9) los tomó por sorpresa. Ellos no esperaban que fuera ahora el profesor quien les hiciera preguntas a ellos, cambiando la posición en este nuevo *sistema de actividad* respecto al anterior (SA2). Inclusive, a partir de este momento, dos estudiantes del equipo se desentienden del diálogo con el profesor.

En SA2, el profesor usa el concepto físico de SR como *recurso* (conceptual) para interpretar el fenómeno físico. Sin embargo, lo usa implícitamente (no les dice a los estudiantes que se va a usar el SR) y de manera poco clara. El SR conlleva una dirección (orientación) y un origen matemático (que puede coincidir o no con el origen fenomenológico). Aquí, el profesor se refiere solamente a la dirección al preguntar "¿para dónde?" (L9). Pero Pedro parece referirse al origen (matemático) del SR, que a su vez lo asocia con el origen fenomenológico, al decir: "desde la parte". Esta falta de entendimiento entre el profesor y Pedro puede deberse a la manera tan espontánea y de manera no explícita, con la que el profesor introduce el concepto del SR. Esto da cuenta de su *trabajo documental* en la manera en que usa el *recurso*. Además, el profesor no da alguna razón del porqué de su pregunta, así como tampoco menciona "*sistema de referencia*" en ningún momento; no obstante que su objetivo es utilizar dicho *recurso* para que Pedro comprenda cómo se establece el signo de la aceleración. Se debe poner atención cuando Pedro dice: "yo lo que logré entender" (L10), ya que Pedro se refiere a la manera en que se llevó a cabo el experimento (SA1); mientras que el profesor pretende llevar a cabo un análisis conceptual y que es lo que guiará el resto de la discusión con Pedro. La discusión continúa:

L11  Prof.: ¿Para dónde tomaste lo positivo? [*Después de que Pedro se queda pensando, el profesor continúa*] las orientaciones, positiva, negativa son arbitrarias, ¿si? Y si tú tomas lo positivo para arriba, no quiere decir que la pelota suba. Y si tomas lo positivo para abajo, no quiere decir que la pelota suba. Y si tomas lo positivo en diagonal, no quiere decir que la pelota suba. La pelota hace su comportamiento y punto [*Pedro dice: "sí"*] ¿No? La otra es la interpretación [*Pedro dice: "Ah, ok."*] ¿Sí?, ¿para dónde tomas tú lo positivo?

L12  Pedro: Pues, este, en el momento en que baja la pelota.

El profesor, usando *recurso* de SR, centra la atención en el observador (al decir: "son arbitrarias") quien es el que elige las convenciones para medir [interpretar] el fenómeno físico (plano inclinado), pero no es explícito que en el experimento, el observador sea el estudiante. Por su parte, Pedro ahora centra su atención en el movimiento del objeto, y comienza a vincularlo con la toma de datos. Sin embargo, no percibe la relación del experimento (movimiento real de la pelota) con el SR. Es aquí cuando el profesor incorpora otro *recurso* de tipo cultural-semiótico: gestos; i.e., un *recurso* utilizado intencionalmente por el profesor para lograr un objetivo (Gueudet & Trouche, 2009) y que a su vez utiliza en un proceso social de creación de significado a través de las acciones del profesor en un tiempo y espacio determinado (Arzarello, Paola, Robutti & Sabena, 2009; Radford 2008). El *uso* de gestos se produce a partir de L13.

L13  Prof.: ¿Para dónde tomas tú lo positivo? No en qué momento, ¿para dónde? [*Pedro le susurra a E1: "¡Ayúdame!" véase Figura 5.1*] Suelto esto [*toma un objeto (mousepad) y lo coloca a una distancia por encima de la mesa; Figura 2a*], ¿qué le va a pasar? [*Pedro dice: "se va a caer"*] [*Prof. Suelta el objeto; véase Figura 5.2b*] ¿Para dónde tomas tú lo positivo? [*Pedro piensa*] Entonces te ayudo. Yo tomo lo positivo para abajo. Entonces, este cuerpo [*refiriéndose al objeto que al mismo tiempo levanta y deja caer*] ¿aumentó su distancia o la disminuyó? [*Pedro no contesta*].

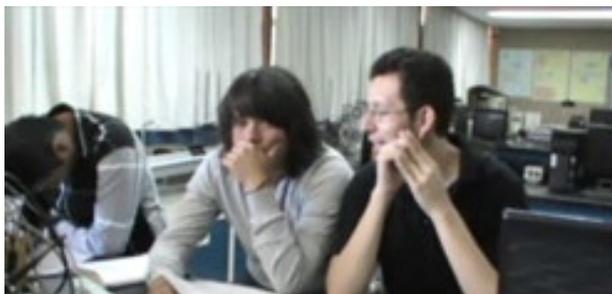

**Figura 5.1. Pedro pide ayuda de un compañero.**

En L13, se observa la dimensión social del momento en que está ocurriendo más allá que una negociación de significados (Radford, 2008) sobre el *recurso* que utiliza el profesor. Pedro se siente nervioso e inseguro al decir: "¡Ayúdame!" (Véase Figura 5.1). Después de que Pedro había respondido lo que él concebía como SR cuando dijo "Lo que yo logré entender" (L10), se encontró en un momento en que no lograba responder lo que el profesor preguntaba. Por lo que la interacción, entre Pedro y el profesor, juega un papel decisivo para reconocer mutuamente los significados institucionales que se están abordando. Esta interacción social que afecta a los individuos es determinante en el proceso de aprendizaje (Radford, 2008) para Pedro y de la manera en que transmite los conocimientos el profesor.

Posteriormente, el *recurso* del gesto es utilizado por el profesor para representar el fenómeno físico y espera que Pedro logre comprender el concepto de SR. En los gestos usados por el profesor, sólo aparece la representación del fenómeno físico (figuras 5.2a y 5.2b) y mediante el lenguaje hace referencia al concepto SR cuando dice: "Yo tomo lo positivo para abajo" (L13). Sin embargo, Pedro parece solamente observar el gesto sin darle significado. Lo que Pedro observa es solamente el movimiento de caída del objeto, la parte visible del *recurso* (Adler, 2000). Cuando el profesor dice: "Yo tomo lo positivo para abajo" en su gesto no hay algún elemento ni de la orientación ni del origen del SR. Por lo que al preguntar si la distancia aumenta o disminuye, el profesor no se da cuenta de que se necesita un origen, respecto del cual aumenta la distancia el objeto. En L14, el profesor modifica –a partir de su reflexión– su *recurso*.

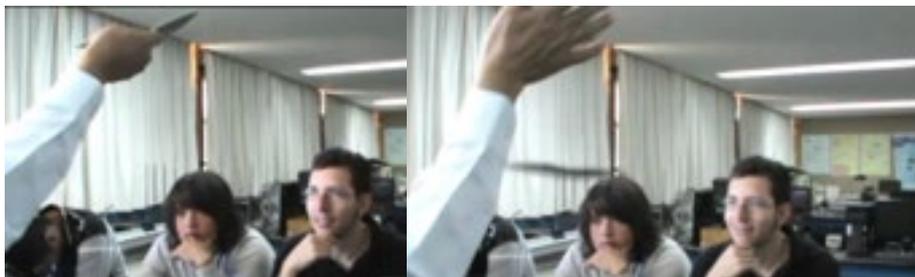

**Figura 5.2.** Gestos como *recursos* utilizados por el profesor en dos momentos (5.2a, izquierda), (5.2b, derecha).

L14   Prof.: ¿A qué distancia está? [*Mientras vuelve a levantar el objeto; Figura 5.2a*] [*Pedro dice: "Este…a equis centímetros"*] Cero, ok. Cero. ¿A qué distancia está? [*Baja el objeto un poco, respecto del punto más alto del que lo había colocado; véase Figura 5.3a*] Aquí estaba el cero [*señala con su dedo dónde se encontraba el objeto en un inicio; véase Figura 5.3b*].

L15   Pedro: Sí, este…podríamos ponerle menos uno o algo así.

L16   Prof.: No, no le puedo poner menos uno. (…) Para abajo es lo positivo [*mientras vuelve a colocar el objeto donde lo había colocado en un inicio y lo deja caer*], yo digo, arbitrariamente. ¿A qué distancia está? [*Coloca el objeto por debajo del punto inicial (cero para el profesor; véase Figura 5.3c*].

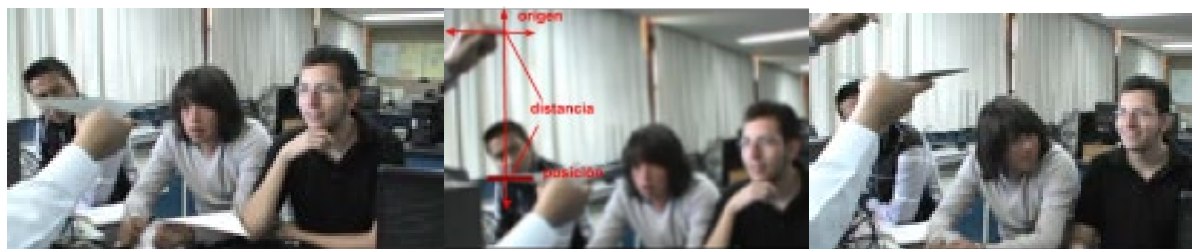

**Figura 5.3.** Gestos como *recursos* utilizados por el profesor en tres momentos (5.3a, izquierda), (5.3b, centro) y (5.3c, derecha).

Mediante la reflexión, referente al *recurso* que utiliza, el profesor (Gueudet & Trouche, 2009) se da cuenta de que necesita un origen respecto del cual se pueda medir la distancia, por lo que incorpora un nuevo gesto (Figura 5.3b) que representa el origen del SR y lo utiliza en dos momentos: cuando dice en L14: "Cero, ok. Cero" y al afirmar "Aquí estaba el cero". Se infiere, la reflexión del profesor sobre el movimiento del objeto, a partir del desarrollo de su discurso usado con los estudiantes, al tratar de explicar el significado del signo de la aceleración. Previamente, no había usado el gesto para representar el origen del SR, pero conforme evoluciona el discurso y la acción del profesor apareció un nuevo gesto (*recurso*). Después de que baja el objeto y pregunta: "¿A qué distancia está?" (Véase Figura 5.3a) es cuando aparece el gesto para hacer referencia al origen del SR (véase Figura 5.3b). Es importante resaltar que el primer momento también indica que el origen del SR (origen matemático) coincide con el origen fenomenológico (inicio del movimiento), pero el profesor no lo dice ni lo hace explícito. El diálogo continúa:

L17   Pedro: Hacia abajo. [*El profesor pregunta nuevamente: "¿a qué distancia está?" haciendo referencia a la Figura 5.3c*] ¡Ay!, no sé.

L18   Prof.: Pues, más o menos calcúlale [*uno de los investigadores de este artículo hace la siguiente intervención: "¿qué distancia hay de su dedo a lo que está sosteniendo?"*].

L19   Pedro: […] Es que ya me puse nervioso, a ver, pongámosle que unos diez centímetros. [*El profesor baja el objeto (lo acerca a la mesa); véase Figura 5.4a*] Como veinticinco. [*El profesor baja aún más el objeto (lo acerca más a la mesa)*] Como treintaicinco, cuarenta, pongámosle cuarenta. [*El profesor sube el objeto y lo coloca por encima de su dedo (origen); véase Figura 5.4b y 5.4c*] Como a diez.

L20   Prof.: Menos diez. [*Al mismo tiempo Pedro dice: "¡ah!, buen punto, menos diez"*].

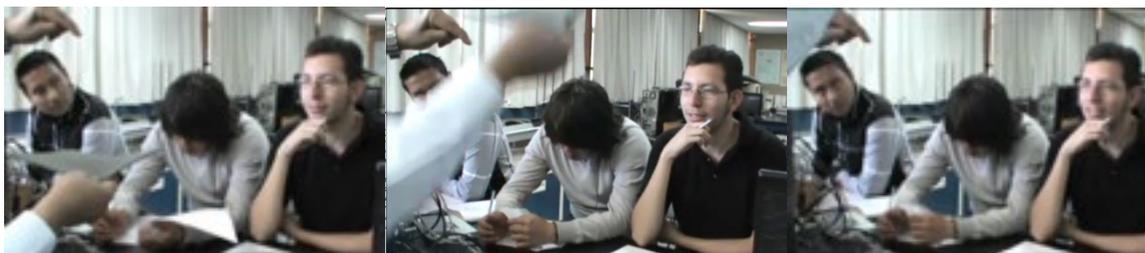

**Figura 5.4.  Gestos como *recursos* utilizados por el profesor en tres momentos para colocar el objeto: a una distancia positiva del origen del SR (5.4a, izquierda) y a una distancia negativa del origen del SR (5.4b, centro) y (5.4c, derecha).**

Pedro no comprende el porqué del gesto del profesor (Figura 5.3b) ni a qué hace referencia. Él sigue haciendo referencia al movimiento que sigue el objeto (fenómeno físico), al decir: "Hacia abajo". Y es hasta que otro integrante de otro equipo le pregunta explícitamente por la distancia que hay del dedo del profesor al objeto, que Pedro logra entender y seguir respondiendo conforme el profesor aleja el objeto de su dedo. Una vez que el profesor determina el SR con su origen hace un nuevo gesto (Figura 5.4b y 5.4c) para tratar de que Pedro se dé cuenta (sin éxito) del papel que juega la orientación; que aquí se representa por el valor negativo [de la aceleración] que tendría el objeto al colocarlo por encima del origen del SR. Es aquí cuando Pedro vuelve a centrarse solamente en la distancia entre el dedo y el objeto (véase Figura 5.4c) sin relacionarla con el concepto de SR.

A continuación, la conclusión de la discusión entre el profesor y Pedro.

### 5.3. Integración de los tres *sistemas de actividad*

El diálogo entre el profesor y Pedro termina cuando el profesor intenta concluir con algunas preguntas a Pedro con el objetivo de determinar si comprendió su explicación previa. En esta última parte (L21-L28) lo que se determina es que Pedro intenta integrar los tres *sistemas de actividad* previos, al relacionarlos sobre el concepto de SR.

L21  Prof.: ¿Sí? [*sobre el valor de menos diez como distancia en la Figura 5.4c*], ¿por qué?

L22  Pedro:  Porque lo positivo lo está tomando hacia abajo.

L23  Prof.: ¿Y eso quiere decir que esto [*refiriéndose al objeto*] se suba? [*Pedro dice: "no"*]. Pues no. ¿Aumentó la distancia o disminuyó?

L24  Pedro:  No, ninguna de las dos.

L25  Prof.: ¿Cómo no?, a ver [*levanta el objeto al punto de inicio (origen) y lo deja caer*] No está en el mismo lugar [*Pedro dice: "Buen punto"*]. ¿Aumenta la distancia o disminuye? [*Vuelve a levantar y dejar caer el objeto*].

L26  Pedro:  Disminuye.

L27  Prof.: ¡Aumenta! ¿No dijimos que para abajo es lo positivo?

L28  Pedro:  ¡Sí es cierto!

En L22 Pedro responde correctamente a partir de la convención utilizada por el profesor (hacia abajo positivo). Sin embargo, cuando en L23 el profesor le pregunta: "¿Aumentó la distancia o disminuyó?", Pedro vuelve a confundirse. Esta confusión, por parte de Pedro, es causada por la pregunta del profesor, pues el hecho de que el profesor haya movido su mano por arriba del origen del SR (Figuras 5.4b y 5.4c), no implica que las distancias entre dos puntos deben disminuir. En física, se pueden obtener valores negativos de las cantidades. Pero, por ejemplo un valor de –10 cm será mayor que –5 cm. Así, aquí el Profesor se sale, de alguna, manera del objetivo central: observar el papel del SR como un conjunto de convenciones para poder obtener medidas de cantidades físicas. De manera que el objetivo común de los tres *sistemas de actividad*: que el análisis del movimiento de objetos, en particular de caída libre, se lleva a cabo a través de las mediciones de cantidades físicas (e.g., posición y distancia) en un SR.

## 6. Conclusiones y reflexiones finales

En este artículo se analizó la manera en que un profesor de física usó gestos como un *recurso* didáctico para atender la duda de un estudiante (E1); y también, como un *artefacto* para representar un concepto: el *sistema de referencia*, el cual a su vez es un *artefacto* en el sentido de Wartofsky (1979). Es un modo de representar la actividad en el cual fue producido, esto es, en el análisis sobre el movimiento de objetos. Esta característica del SR permitió que este concepto pudiera ser el elemento principal sobre el cual se vincularon –por parte del profesor– tres *sistemas de actividad*.

Así, respecto de la primera pregunta de investigación, el profesor representó a través del *uso* de gestos: el origen de un SR y cómo varían las unidades de la distancia cuando un objeto se mueve en ese SR. Así se pudo observar y resaltar la componente didáctica de los gestos. Respecto de la segunda pregunta de investigación, el análisis del *uso* de gestos en el profesor permitió dar cuenta del conocimiento que él tiene sobre el concepto de SR. Se observó que identifica los elementos esenciales del SR: la necesidad de un origen matemático, de un observador quien lleva a cabo las mediciones, y la posibilidad de cambiar las orientaciones de medición (positiva y negativo) de un SR. Así, se resalta la componente epistemológica y cognitiva de los gestos.

Es importante señalar que, de acuerdo al encuadre teórico utilizado (§ 2), no se analizó la estructura cognitiva [conocimiento] del profesor a través de las invariantes operatorias del *uso* de gestos; sino que considerando a los gestos como un medio cognitivo y un modo de representación, se pudo así inferir el conocimiento del profesor. Esto, debido a la componente semiótica de los gestos como portadores de significados. Así, este estudio sirve como argumento para manifestar el hecho de que los gestos deben ser considerados como parte importante del *sistema de recursos* de los profesores.

De manera que futuras investigaciones estarían encaminadas en analizar la manera en que los gestos se integran al *sistema de recursos* del profesor y en cómo se desarrollan a través de la *génesis documental* para determinar, por ejemplo, el papel de los gestos en el desarrollo de una lección. Para esto, sería necesario incorporar tanto en el diseño de las investigaciones como en el análisis: entrevistas a los profesores y la manera en que los gestos se integran con otros *recursos* y en variedad de situaciones.

## Bibliografía

**Ulises Alfonso Salinas-Hernández**

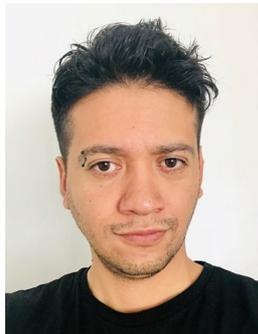Doctor en Ciencias con Especialidad en Matemática Educativa por el Cinvestav-IPN Zacatenco, México. Profesor de asignatura en la Escuela Nacional Colegio de Ciencias y Humanidades de la UNAM, México. Interesado en analizar la enseñanza-aprendizaje de las matemáticas y de la ciencia como una práctica social, histórica y culturalmente mediada. Profesor visitante en *L'École des sciences de l'éducation de l'Université Laurentienne* en 2015, en donde trabajó con el profesor Luis Radford. Realizó una estancia doctoral en *The French Institute of Education, École Normale Supérieure de Lyon, France* en 2018, bajo la supervisión del profesor Luc Trouche. Ha participado en diversos congresos internacionales en el ámbito de la matemática educativa. E-mail: asalinas@cinvestav.mx

**Luc Trouche**

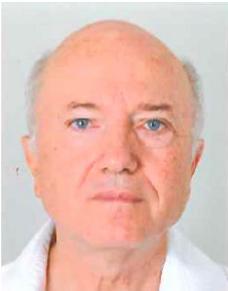Luc Trouche is professor of didactics of mathematics in the French Institute of Education, École Normale Supérieure de Lyon, France. Interested in studying the role of tools in the learning of mathematics (Monaghan, Trouche & Borwein, 2016), he has contributed to develop the notion of instrumental orchestration (Trouche & Drijvers, 2014) for modeling the management, by a teacher, of available tools in the classroom. He focuses now on resource use/design and teacher professional development in the time of digitalization. This has led him to contribute to develop the documentational approach to didactics (Gueudet, Pepin & Trouche, 2012). In this perspective, the notion of teacher resource system appears crucial in order to understand teacher (developing) knowledge and the coherence of his/her activity. Studying the interactions between individual and collective teachers' resource systems gives means for understanding the dynamics of these collectives, and for rethinking the way of supporting teacher development at a time of the 'metamorphosis' of teaching environments. E-mail: luc.trouche@ens-lyon.fr